\documentclass[12pt]{amsart}
\usepackage{amsmath,amsthm}

\def\pmatrix{\left(\begin{matrix}}
\def\endpmatrix{\end{matrix}\right)}

\def\P{{\mathbb P}}
\def\Z{{\mathbb Z}}
\def\C{{\mathbb C}}

\def\tn{{\theta_{\rm null}}}

\def\de{\delta}

\def\p{\partial}

\def\t{\theta}

\def\T{\Theta}
\def\e{\varepsilon}

\def\O{{\mathcal O}}
\def\A{{\mathcal A}}

\def\J{{\mathcal J}}
\def\H{{\mathcal H}}
\def\D{{\mathcal D}}
\def\B{{\mathcal B}}
\def\m{{\mathfrak m}}

\def\tch#1#2{{\left[\begin{matrix}{#1}\\ {#2}\end{matrix}\right]}}
\def\tt#1#2{{\t\tch#1#2}}
\def\Sp#1{\operatorname{Sp}(#1,\Z)}
\def\Sing{\operatorname{Sing}}
\def\tn{\theta_{\rm null}}

\theoremstyle{plain}
\newtheorem{thm}{Theorem}%[section]
\newtheorem{lm}[thm]{Lemma}

\newtheorem{cor}[thm]{Corollary}

\theoremstyle{definition}
\newtheorem{df}[thm]{Definition}

\newtheorem{dsc}[thm]{}

\title{Jacobians with a vanishing theta-null in genus 4}
\author{Samuel Grushevsky\and Riccardo Salvati Manni}
\address{Mathematics Department, Princeton University, Fine Hall,
Washington Road, Princeton, NJ 08544, USA}
\email{sam@math.princeton.edu}
\address{Dipartimento di Matematica, Universit\`a ``La Sapienza'',
Piazzale A. Moro 2, Roma, I 00185, Italy}
\email{salvati@mat.uniroma1.it}
\begin{document}
\begin{abstract}
In this paper we prove a conjecture of Hershel Farkas \cite{farkas}
that if a 4-dimensional principally polarized abelian variety has a
vanishing theta-null, and the  hessian of the theta function at the
corresponding point of order two is degenerate, the abelian
variety is a Jacobian.

We also discuss possible generalizations to higher genera, and an
interpretation of this condition as an infinitesimal version of
Andreotti and Mayer's local characterization of Jacobians by the
dimension of the singular locus of the theta divisor.
\end{abstract}
\maketitle

\section*{Introduction}
The study of the geometry of the theta divisor of Jacobians of
algebraic curves is a very classical subject going back at least to
Riemann's celebrated theta-singularity theorem itself, later
strengthened by Kempf \cite{kempf}. The geometry of the canonical
curve was also shown to be related to the geometry of the theta
divisor --- in particular the celebrated Green's conjecture
\cite{green} says that the tangent cones to the Jacobian theta
divisor at its double points span the ideal of quadrics containing
the canonical curve.

These tangent cone quadrics in fact have rank 4; it can also be
shown that the singular locus of the theta divisor $\Sing\T$ of a
Jacobian of a curve of genus $g$ has dimension $\ge g-4$ ($g-3$ for
hyperelliptic curves). It was thus asked whether this is a
characteristic property for the locus of Jacobians $\J_g$ within the
moduli space of principally polarized abelian varieties.

This study was undertaken by Andreotti and Mayer \cite{am}. Let
$\A_g$ be the moduli space of (complex) principally polarized
abelian varieties --- ppavs for short. Denote by $N_k\subset\A_g$
the locus of ppavs for which $\dim\Sing\T\ge k$. Andreotti and Mayer
showed that $\J_g$ is an irreducible component of $N_{g-4}$. The
situation was further studied by Beauville \cite{gen4am} and Debarre
\cite{debarredecomposes}. It was shown that
for $g\ge 4$ the locus $N_{g-4}$ is reducible (and thus not equal to
$\J_g$); however, conjecturally all components of $N_{g-4}$ other
than $\J_g$ are contained in the theta-null divisor $\tn$ (the zero
locus of the product of all theta constants with half-integral
characteristics; alternatively, the locus of those ppavs for which
the theta divisor has a singularity at a point of order two), which
is a component of $N_0$.

Thus it is natural to try to study the intersection of $\J_g$ with
the other components of $N_{g-4}$, or at least with $\tn\subset
N_0$. This is the object of this paper: we study the infinitesimal
version of $N_{g-4}$, prove H. Farkas' conjecture \cite{farkas}
describing the situation for $g=4$, and discuss the possible
situation in general.

In a recent paper \cite{diffeq}, we proved some identities between
the simplest types of theta series with harmonic polynomial
coefficients, which are generalizations of the classical Jacobi's
derivative formula. Proving these involves the evaluation at zero of
the derivatives of first and second order of classical theta
functions, i.e. some local infinitesimal properties of the theta
divisor. The methods developed there are applicable to our problem
and provide a natural generalization of Farkas' question as follows.

Consider the following stratification of $\tn$: let
$\tn^{h}\subset\tn$ be the subset where the tangent cone to the
theta divisor at the corresponding singular point of order two has
rank $\le h$, for $h=0,1,\dots,g$. We have equations for the (level
covers of) loci $\tn^{h}$, and believe that a further investigation
of the relation between $\J_g\cap\tn$ and $\tn^{3}$ would be of
interest. Similarly it seems promising to study the locus $\tn^2$
locally near the locus of reducible ppavs (which is contained in
it). In this paper we treat the $g=4$ case completely; higher
dimensions will be addressed in a separate paper.

\section{Notations and definitions}
\begin{df}
Let $\H_g$ denote the {\it Siegel upper half-space}, i.e. the set of
symmetric complex $g\times g$ matrices $\tau$ with positive definite
imaginary part. Each such $\tau$ defines a complex principally
polarized abelian variety (ppav for short)
$A_\tau:=\C^g/\Z^g+\tau\Z^g$. If $M=\pmatrix a&b\\
c&d\endpmatrix\in\Sp{2g}$ is a symplectic matrix in a $g\times g$
block form, then its action on $\tau\in\H_g$ is defined by
$M\circ\tau:=(a\tau+b)(c\tau+d)^{-1}$, and the moduli space of ppavs
is the quotient $\A_g=\H_g/\Sp{2g}$.

A period matrix $\tau$ is called {\it reducible} if there exists
$M\in\Sp{2g}$ such that
$$
  M\cdot\tau=\pmatrix \tau_1&0\\
  0&\tau_2\endpmatrix,\quad\tau_i\in\H_{g_i},\ g_1+g_2=g;
$$
otherwise we say that $\tau$ is irreducible.
\end{df}
\begin{df}
For $\e,\de\in (\Z/2\Z)^g$, thought of as vectors of zeros and ones,
$\tau\in\H_g$ and $z\in \C^g$, the {\it theta function with
characteristic $[\e,\de]$} is
$$
  \tt\e\de(\tau,z):=\sum\limits_{m\in\Z^g} \exp \pi i \left[
  ^t(m+\frac{\e}{2})\tau(m+\frac{\e}{2})+2\ ^t(m+\frac{\e}{2})(z+
  \frac{\de}{2})\right].
$$

A {\it characteristic} $[\e,\de]$ is called {\it even} or {\it odd}
depending on whether $\tt\e\de(\tau,z)$ is even or odd as a function
of $z$, which corresponds to the scalar product
$\e\cdot\de\in\Z/2\Z$ being zero or one, respectively. A {\it theta
constant} is the evaluation at $z=0$ of a theta function. All odd
theta constants of course vanish identically in $\tau$.

Observe that
$$
  \tt 0 0\left(\tau,z+\tau\frac\e2+\frac\de2\right)= \exp \pi i
  \left(-\frac{^t\e}{2}\,\tau\frac{\e}{2}\,\,-\frac{^t\e}{2}\,(z+
  \frac{\de}{2})\right)\tt\e\de(\tau,z),
$$
i.e. theta functions with characteristics are, up to some easy
factor, the Riemann's theta function (the one with characteristic
$[0,0]$) shifted by points of order two.
\end{df}

\smallskip
Let $\rho:{\rm GL}(g,\C)\to \operatorname{End} V$ be an irreducible
rational representation with the highest weight
$(k_1,k_2,\dots,k_g)$, $k_1\geq k_2 \geq\dots\geq k_g$; then we call
$k_g$ the weight of $\rho$. A representation $\rho_0$ is called
reduced if its weight is equal to zero. Let us fix an integer $r$;
we are interested in pairs $\rho=(\rho_0,r)$, with $\rho_0$ reduced.
We call $r$ the weight of $\rho$ and use the notation
$$
  \rho(A)=\rho_0(A)\det A^{r/2}\ .
$$
\begin{df} A map $f:\H_g\to V$ is called a {\it $\rho$-valued modular
form} with respect to a finite index subgroup $\Gamma\subset\Sp{2g}$
if
$$
  f(\sigma\circ\tau)=v(\sigma)\rho(c\tau+d)f(\tau)\qquad\forall
  \tau\in\H_g,\forall\sigma\in\Gamma,
$$
and if additionally $f$ is holomorphic at all cusps of
$\H_g/\Gamma$.

If $\rho(\sigma)=\det(c\tau+d)^k$, then we call this a {\it scalar
modular form} of weight $k$.
\end{df}

For a finite index subgroup $\Gamma\subset\Sp{2g}$ a multiplier
system of weight $r/2$ is a map $v:\Gamma\to \C^*$, such that the
map
$$
  \sigma\mapsto v(\sigma)\det(c\tau+d)^{r/2}
$$
satisfies the cocycle condition for every $\sigma\in\Gamma$ and
$\tau\in\H_g$ (note that the function $\det(c\tau+d)$ possesses a
square root). Clearly a multiplier system of integral weight is a
character.

\begin{df} A map $f:\H_g\to V$ is called a $\rho$- or {\it $V$-valued modular
form}, or simply a {\it vector-valued modular form}, if the choice
of $\rho$ is clear, {\it with multiplier $v$}, with respect to a
finite index subgroup $\Gamma\subset\Sp{2g}$ if
\begin{equation}\label{transform}
  f(\sigma\circ\tau)=v(\sigma)\rho(c\tau+d)f(\tau)\qquad\forall
  \tau\in\H_g,\forall\sigma\in\Gamma,
\end{equation}
and if additionally $f$ is holomorphic at all cusps of
$\H_g/\Gamma$.
\end{df}

\begin{df}[Theta constants are modular forms]
We define the {\it level} subgroups of the symplectic group to be
$$
  \Gamma_g(n):=\left\lbrace M=\pmatrix a&b\\ c&d\endpmatrix
  \in\Sp{2g}\, |\, M\equiv\pmatrix 1&0\\
  0&1\endpmatrix\ {\rm mod}\ n\right\rbrace
$$
$$
  \Gamma_g(n,2n):=\left\lbrace M\in\Gamma_g(n)\, |\, {\rm
  diag}(a^tb)\equiv{\rm diag} (c^td)\equiv0\ {\rm mod}\
  2n\right\rbrace.
$$

Under the action of $M\in\Sp{2g}$ the theta functions transform as
follows:
$$
  \theta \bmatrix M\pmatrix \e\\ \de\endpmatrix\endbmatrix
  (M\cdot\tau,\,^{t}(c\tau+d)^{-1}z)\qquad\qquad\qquad
$$
$$
  \qquad\qquad\qquad=\phi(\e,\,\de,\,M,\,
  \tau,\,z)\det(c\tau+d)^{\frac{1}{2}}\theta\bmatrix \e\cr \de
  \endbmatrix(\tau,\,z),
$$
where
$$
  M\pmatrix \e\cr \de\endpmatrix :=\pmatrix d&-c\cr
  -b&a\endpmatrix\pmatrix \e\cr \de\endpmatrix+ \pmatrix {\rm
  diag}(c \,^t d)\cr {\rm diag}(a\,^t b)\endpmatrix
$$
taken modulo 2, and $\phi(\e,\,\de,\,M,\,\tau,\,z)$ is some
complicated explicit function. For more details, we refer to
\cite{igbook} and \cite{farkasrauch}.

Thus theta constants with characteristics are (scalar) modular forms
of weight $1/2$ with respect to $\Gamma_g(4,8)$, i.e. we have
$$
  \tt\e\de(M\circ\tau,0)=\det(c\tau+d)^{1/2}\tt\e\de(\tau,0)
  \qquad \forall M\in\Gamma_g(4,8).
$$
\end{df}
\begin{df} We call the {\it theta-null divisor} $\tn\subset\A_g$ the
zero locus of the product of all even theta constants. We define a
stratification of $\tn$ as follows. For $h=0,\ldots,g$ we let
$$
  \tn^h=\left\lbrace\tau\in\H_g: \exists[\e,\de]\ {\rm even}\ ,
  \tt\e\de(\tau)=0;\ \operatorname{rk}\left.
  \frac{\partial^2\tt\e\de(\tau,z)}{\partial z_i\partial
  z_j}\right|_{z=0}\le h\right\rbrace,
$$
i.e. the locus of points on $\tn$ where the rank of the tangent cone
to the theta divisor at the corresponding point
$\frac{\tau\e+\de}{2}$ of order two is at most $h$.
\end{df}

By the above transformation formulae, we see that $\tn$ and
$\tn^{h}$ are well-defined on $\A_g$ and not only on the {\it level
moduli spaces} $\A_g(4,8):=\H_g/\Gamma_g(4,8)$. Since the theta
constant with characteristic is up to a non-zero factor the value of
Riemann's theta function at the corresponding point of order two,
$\tn$ can also be described as the locus of ppavs for which $\T$ has
a singularity at a point of order two. Note that $\t(\tau,z)$ is an
even function of $z$ (the theta divisor is symmetric under the
involution $\pm 1$), and thus the locus $\tn$ indeed turns out to be
a divisor).

\begin{dsc}[The ring of modular forms]
We recall that the theta constants define an embedding
$$
  \begin{aligned}
  Th:\A_g(4,8)&\rightarrow\P^{2^{g-1}(2^g+1)-1}&\\
  \tau&\mapsto\left\lbrace\tt\e\de(\tau)\right\rbrace_{[\e,\de]\ {\rm
  even}}& \end{aligned}.
$$
This map extends to the Satake compactification
${\overline{\A_g(4,8)}}$. Hence the ring of scalar modular forms for
$\Gamma(4,8)$ is the integral closure of the ring
$\C\left[\tt\e\de\right]$. The ideal of algebraic equations defining
$Th({\overline{\A_g(4,8)}})\subset \P^{2^{g-1}(2^g+1)-1}$ is known
completely only for $g\le 2$ (and almost known for $g=3$, see
\cite{vgvdg}).
\end{dsc}

\begin{dsc}
Since theta functions satisfy the heat equation
$$
 \frac{\partial^2\tt\e\de(\tau,z)}{\partial z_i\partial z_j}
 =\pi i(1+\delta_{i,j})\frac{\partial\tt\e\de(\tau,z)}{\partial\tau_{ij}},
$$
(where $\delta_{i,j}$ is Kr\"onecker's delta), the Hessian of the
theta funtions with respect to $z$ can be rewritten as the first
derivatives with respect to $\tau_{ij}$. Hence if a point
$x=\tau\frac\e2+\frac\de2$ of order two  is a singular point in the
theta divisor, which is simply to say $\tt 00(\tau,
x)=0=\tt\e\de(\tau,0)$ (the first derivatives at zero of an even
function are all zero), the rank of the quadric defining the tangent
cone at $x$ is the rank of the matrix obtained by applying the
$g\times g$-matrix-valued differential operator
$$
  \D:= \left(\begin{array}{rrrr}
  \,\frac{\p}{\p\tau_{11}}&\frac{1}{2}\frac{\p}
  {\p\tau_{12}}&\dots&\frac{1}{2}\frac{\p}{\p\tau_{1 g}}\\
  \frac{1}{2}\frac{\p}{\p \tau_{21}}&\frac{\p}{\p
  \tau_{22}}&\dots&\frac{1}{2}\frac{\p}{\p\tau_{2 g}}\\
  \dots&\dots&\dots&\dots\\
  \frac{1}{2}\frac{\p}{\p \tau_{g 1}}& \dots&\dots& \,\
  \frac{\p}{\p\tau_{g g}}\end{array}\right)
$$
to  $\tt\e\de(\tau,0)$.
\end{dsc}

\section{Equations for $\tn^h$}
The locus $\tn^{h}$ is given by the conditions
$$
 \lbrace  \exists\,\, [\e,\de]\ {\rm even;}\ 0=\tt\e\de
 (\tau);\ {\rm rk}\, \D\tt\e\de(\tau)\le h\rbrace.
$$
We can get equations for $\tn^h$ by setting all $(h+1)\times (h+1)$
minors of $\D\tt\e\de(\tau)$ equal to zero, but these minors are not
modular forms: the derivative of a section of a bundle is only a
section of that bundle when restricted to the zero set of the
section, i.e. $\D\tt\e\de(\tau)$ is not a modular form, but is
modular when restricted to the locus $\tt\e\de(\tau)=0$. This
condition is not invariant under $\Sp{2g}/\Gamma_g(4,8)$, and thus
for technical reasons we will work on $\A_g(4,8)$. However, the
locus $\tn^h$ that we are describing is invariant under the action
of $\Sp{2g}$, and we will thus be able to easily descend from
$\A_g(4,8)$ to $\A_g$ by symmetrizing.

The divisor $\tn\subset \A_g(4,8)$ is reducible. Its irreducible
components are the divisors of individual theta constants with
characteristics --- cf. \cite{fr} page 88 for $g\geq 3$ and by
inspection in the remaining  two cases . These components are all
conjugate under the action of $\Sp{2g}$, and thus for our purposes
we can restrict to one component. Without loss of generality we can
take this to be $\t_0:=\lbrace \tt 0
0(\tau)=0\rbrace\subset\A_g(4,8),$ and consider its stratification,
letting $\t_0^h$ be those $\tau\in\t_0$ for which the rank of the
tangent cone, i.e. of the hessian of the theta function, at zero is
at most $h$.

Following the ideas of \cite{diffeq}, we observe that for any even
characteristic $[\e,\,\de]\neq [0,\,0] $ the expression
$$
  \tt\e\de(\tau)^2\D(\tt 0 0/\tt\e\de)(\tau)
$$
$$
  =(1+\delta_{i,j})\left[\tt\e\de(\tau)\frac{\partial\tt00(\tau)}
  {\partial\tau_{ij}}-\tt 0 0(\tau)\frac{\partial\tt\e\de(\tau)}
  {\partial\tau_{ij}}\right]
$$
is a vector-valued modular form of weight $1$ with
$\rho_0=(2,0,\dots,0)$ with respect to $\Gamma_g(4,8)$.

We denote by $\B\left(\tch00,\tch\e\de\right)^h(\tau)$ the
$\pmatrix g\\ h\endpmatrix \times\pmatrix g\\
h\endpmatrix\ $ symmetric matrix obtained by taking in lexicographic
order all the $h\times h$ of the above matrix --- it is a vector
valued modular form of weight $h$ with $\rho_0=(2,\dots,2,
0\dots,0)$, with respect to $\Gamma_g(4,8)$.
\begin{thm}\label{th}
The locus $\t_0^h$ is set theoretically defined by
$$
  \tt 0 0(\tau,0)=0=\B\left(\tch00,\tch\e\de\right)^{h+1}(\tau)
  \qquad\forall [\e,\de]\ne[0,0]\ {\rm even}.
$$
\end{thm}
\begin{proof}
One implication is trivial. Viceversa, let us assume that all
equations are satisfied; there always exists an even characteristic
$[\e,\de]$ such that $\tt\e\de(\tau)\neq 0$, thus
$\B\left(\tch00,\tch\e\de\right)^{h+1}(\tau)=0$ implies
$\tau\in\t_0^h$.
\end{proof}

\section{Infinitesimal Andreotti-Mayer condition }
We would like to understand the loci $\tn^h$ and $\tn\cap\J_g$.
Suppose we have a Jacobian with a vanishing theta-null, i.e.
$X\in\J_g\cap\tn$, and $x\in X[2]$ is the corresponding point of
order two on the theta divisor. It is a consequence of Kempf's
singularity theorem (see \cite{kempf}, \cite{ac}) that the rank of
the hessian of the theta function (i.e. of the tangent cone) at $x$
is equal to three. By Green's conjecture \cite{green} the quadrics
defining the canonical curve are obtained exactly as Hessians of the
theta function with respect to $z$ at singular points of the theta
divisor. In general they have rank at most 4, but it is easy to see
that at singular points of order two the rank actually drops to
three --- see \cite{ac} and also \cite{farkas}. Thus $X$ belongs to
$\tn^3$, so we see that $(\J_g\cap\tn)\subset\tn^3$.

It is natural to view the Hessian rank condition $\tn^h$ as the
infinitesimal version of the Andreotti-Mayer condition, indicating
that $\Sing\,\T$ has a $(g-4)$-dimensional 2-jet at $x$.

The codimension of the space of $g\times g$ symmetric matrices of
rank at most 3 within $\H_g$ is equal to $(g-1)(g-2)/2$ (the
dimension is $3g-3$: choosing the first three rows, and thus also
the first three columns, generically determines everything), and
thus for dimension reasons it is natural to expect that locally the
condition of the rank being 3 characterizes Jacobians with a
vanishing theta-null. We believe that the infinitesimal
Andreotti-Mayer condition at points of order two should help to
characterize Jacobians with a vanishing theta-null locally. We thank
Enrico Arbarello, Hershel Farkas and Edoardo Sernesi for stimulating
discussions on this topic.

To study this link , one could try to ``integrate'' the local
condition, to show that ${\rm Sing}\,\Theta$, indeed, has dimension
$g-4$, but this seems hard. By using the ideas of \cite{am} (see
also \cite{ac}), we can get for $\tau\in\tn^3$ some equations
involving the  second derivative of the theta functions with respect
to the variable $\tau$, but we could not give a good dimensional
estimate of the conditions thus obtained.

Another potential approach to this problem would be to use
theorem 6 from \cite{diffeq} to express the $\B$ in terms of
Jacobian determinants of odd theta functions, with syzygetic
characteristics. If one could then prove the appropriate
generalization of Jacobi's derivative formula, conjectured for all
genera and proven for genus 4 in \cite{nonazyg}, this expression in
terms of Jacobian determinants could then be rewritten as an
algebraic expression in terms of theta constants, and then could
perhaps be compared to Schottky-Jung \cite{farkasrauch} equations
for theta constants, or could be used to give a conjectural local
algebraic solution to the Schottky problem.

One could also try to study the infinitesimal Andreotti-Mayer
condition at singular points of the theta divisor that are {\it not}
points of order two, and see whether the rank of the Hessian there
can be used to locally characterize the Jacobian locus. In genus 4
there is no problem, since $N_0=\tn\cup\J_4$, but the situation in
higher dimensions seems very complicated, as it is not clear how to
interpret such a condition in terms of modular forms.

It is known that $N_{g-2}$ is the locus of reducible ppavs
\cite{el}. It is then clear that ppavs with reducible theta divisor
are in $\tn^2$ --- in this case the tangent cone is a quadric that
is the union of two hyperplanes. We are naturally led to consider
the link between $N_{g-2}$ and $\tn^2$. This will be done in  a
forthcoming paper.

\section{The theorem in genus 4}
Our main result of this section is the proof of the following
conjecture of H. Farkas \cite{farkas}:
\begin{thm}\label{mainthm}
If for a 4-dimensional ppav $\tau\in\A_4$ some theta constant and
its hessian are both equal to zero, this ppav is a Jacobian, i.e.
$$
  \tn^3=\J_4\cap\tn.
$$
\end{thm}
\begin{proof}
The fact that this vanishing holds for Jacobians (i.e. the
implication $\Leftarrow$) is discussed in the previous section.

Since $\Sp{4}/\Gamma_4(4,8)$ acts transitively on the set of theta
constants with characteristics, it is enough to take
$[\e,\de]=[0,0]$ above, and show (here $\t$ denotes the theta
function with zero characteristics) that if
$\t(\tau)=\det\D\t(\tau)=0$, then $\tau\in\J_4\cap\tn$.

We denote by $\J_4(4,8)\subset\A_4(4,8)$ the Jacobian locus; also
denote $J:=Th(\J_4(4,8))\subset A:=Th(\A_4(4,8))$, and denote the
chosen component of the theta-null by $T:=A\cap \lbrace
\theta(\tau)=0\rbrace$. Let us denote the locus we are interested in
by $S:=T\cap\lbrace\det\D\theta(\tau)=0\rbrace$.

The main theorem is then the statement that $S=J\cap T$, of which we
already know the inclusion $S\supset J\cap T$ from the previous
discussion.
\begin{lm}
$\dim S=\dim(J\cap T)=8$.
\end{lm}
\begin{proof}
It is known that no theta constant vanishes identically on the
Jacobian locus. Thus $J\cap T\not\subseteq J$, and is given by one
equation, so $\dim (J\cap T)=\dim J-1=8$. On the other hand, the
locus $\lbrace\det\D\theta(\tau)=0\rbrace\subset \H_4$ does not
contain $\lbrace\theta(\tau)=0\rbrace\subset \H_4$, since they are
both of codimension 1 in $\H_4$, and the first locus is not
invariant under $\Gamma_4(4,8)$, while the second is. Thus we have
$S\not\subseteq T$, and since locally in $\H_g$ it is given by one
extra equation, we get $\dim S=\dim T-1=\dim A-2=8$.
\end{proof}
By the discussion in the previous section we know that $S\supset
J\cap T$, and to prove that $S=J\cap T$ it is enough to show that
the degrees of the two sets, as of $8$-dimensional subvarieties of
the projective space, are equal. We will now compute these degrees.
\begin{lm}
The set-theoretic degree $\deg(J\cap T)=8\deg A$.
Scheme-theoretically, the degree of this intersection is $16\deg A$.
\end{lm}
\begin{proof}
Since $\lbrace\theta(\tau)=0\rbrace$ is a hyperplane in $\P^{135}$,
we have $\deg T=\deg A$, and $\deg(J\cap T)=\deg J$. Recalling that
$A$ and $J$ are both irreducible varieties, and $J\subset A$ is
given by one equation of degree 16 in theta constants, i.e. by a
polynomial of degree 16 \cite{vgvdg} in the coordinates of
$\P^{135}$, we see that $\deg J=16\deg A$ scheme-theoretically.

However, the equation for $J\subset A$ in the notations of
\cite{vgvdg} is
$$
  (r_1+r_2+r_3)(r_1-r_2+r_3)(r_1+r_2-r_3)(r_1-r_2-r_3)=0,
$$
where each $r_i$ is a monomial of degree 4 in theta
constants.\footnote{There are many possible choices for $r_i$, and
thus many resulting forms of the equation, which are all conjugate
under the action of $\Sp{4}/\Gamma_4(4,8)$. What happens is that
$A\subset\P^{135}$ is itself given by a large number of equations
(explicitly unknown to this date), and when we intersect $A$ with
{\it any} equation of the above form, the intersection is always the
same, and in particular invariant under $\Sp{4}/\Gamma_4(4,8)$.} If
a theta constant vanishes, i.e. if we are on the theta-null divisor,
it means that one of the products $r_i$ vanishes, so without loss of
generality let's say $r_3=0$. In this case the above equation
becomes $(r_1-r_2)^2(r_1+r_2)^2=0$, so that there is multiplicity at
least two. To finish the proof of the lemma, we need to show that
the multiplicity is exactly two, for which we will use the argument
similar to the one used in \cite{irr} to prove local irreducibility
of Schottky's divisor.

Indeed, take a diagonal period matrix $\tau_0\in\H_1 ^{\times 4}$
with diagonal entries $\omega_1, \omega_2,\omega_3,\omega_4$, and
denote by $\O$ the analytic local ring of $\H_4$ at $\tau_0$. We
shall use $\tau_{ii}-\omega_i$ for $1\leq i\leq 4$, and $2\pi {\sqrt
-1} \tau_{ij}$ for $1\leq i< j\leq 4$ as the generators of the
maximal ideal $\m\subset\O$. If we arrange the $2\pi {\sqrt -1}
\tau_{ij}$ in the order $(12), (34), (13), (24), (14), (23)$ and
call them $x_1, x_2,\dots, x_6$, then the expression
$(r_1+r_2+r_3)(r_1-r_2+r_3)(r_1+r_2-r_3)(r_1-r_2-r_3)$ in
$\m^8/\m^9$ is equal to
$$
  2^{16}\delta(\omega_1)\delta(\omega_2)\delta(\omega_3)\delta(\omega_4)
  P(X),
$$
where $\delta(\omega)$ is the unique cusp form of weight $12$ for
$\operatorname{SL}(2,\Z)$, suitably normalized, cf. \cite{irr}  and
$$
  P(X)=(x_1x_2-x_3x_4)^2 (x_5x_6)^2 - 2(x_1x_2+x_3x_4)\prod_{i=1}^6
  x_i +\left(\prod_{i=1}^4 x_i\right)^2.
$$
For simplicity choose, for this lemma only, the vanishing theta
constant to be $\tt\e\de(\tau)$ with $\e=\de=(1\,1\,0\,0)$, instead
of the zero characteristic --- this theta constant in $\m/\m^2$ has
local expression $x_1$. Hence, by substitution, we have the
expression $(x_3x_4x_5x_6)^2$ for $P(X)$, so the multiplicity of the
intersection of $T$ and $J$ (as given by the Schottky relation) is
exactly two. As a consequence, the set-theoretic degree is half of
the scheme-theoretic degree.
\end{proof}

To compute $\deg S$, we use results from the previous sections. Our
original problem in dealing with $S$ is that $\det\D(\theta(\tau))$
is not a modular form, and thus its zero locus is not well-defined
on $\A_g(4,8)$. However, we can apply theorem \ref{th} for $g=4,
h=3$.
\begin{lm}\label{lm}
Define a function $F$ on $\H_g$ by
$$
  F(\tau):=\left(\tt\e\de(\tau)\right)^{2g}\det\D
  \left(\frac{\tt 0 0(\tau)}{\tt\e\de(\tau)}\right).
$$
Then for any even $[\e,\de]$, scheme-theoretically we have
$$
  T\cap\lbrace F(\tau)=0\rbrace =T\cap\left\lbrace0=\left(\tt\e\de(\tau)\right)^{g}
  \det\D\theta(\tau)\right\rbrace
$$
\end{lm}
\begin{proof}
This is obtained by writing out the derivatives involved in $F$
and using $\tt00(\tau)=0$ on $T$.
\end{proof}

As discussed in \cite{diffeq}, $F(\tau)$ is a scalar modular form of
weight $g+2$ with respect to $\Gamma_g(4,8)$ (recall that theta
constants have weight $1/2$). We can now compute the degree of $S$
and thus finish the proof of the theorem.
\begin{lm}
The set-theoretic degree $\deg S=8\deg A\ (=\deg( J\cap T))$.
\end{lm}
\begin{proof}
For $g=4$ the form $F$ is of weight 6, and thus it follows that it
is a section of the 12'th power of the bundle of theta constants,
the degree of the zero locus of $F$ in $\P^{135}$ is equal to 12.
Thus we have $\deg(T\cap\lbrace F(\tau)=0\rbrace)=12\deg T=12\deg
A$.

To understand the locus $T\cap\lbrace F(\tau)=0\rbrace$, note that
on the right-hand-side of lemma \ref{lm} we have the union of two
loci: $S$, which is exactly $T\cap\lbrace\det\D\theta(\tau)
=0\rbrace$, and of $T\cap\lbrace \tt\e\de(\tau)=0\rbrace $, with
multiplicity $g=4$. Of course the latter is the intersection of $A$
with two hyperplanes, so its degree is equal to $\deg A$. Thus
comparing the scheme-theoretic degrees on both sides of the equality
in the previous lemma, we get the scheme-theoretic degree $\deg
S=12\deg A-4\deg A=8\deg A$. Since $S\supset (J\cap T)$ of the same
dimension, which already has degree $8\deg A$, there is no
multiplicity, and the set-theoretic degree of $S$ is the same as
scheme-theoretic.
\end{proof}
Since we have shown that $\deg S=\deg J\cap T$ as of sets, and we
know $S\supset (J\cap T)$, the theorem --- the statement $S=J\cap T$
--- finally follows.
\end{proof}

As an immediate consequence we have the following
\begin{cor}
If for $\tau\in\A_4$ some theta constants and all minors of order 3
of its hessian are equal to zero, this ppav has reducible theta
divisor.
\end{cor}
\begin{proof}
By the theorem we know that $\tau\in\overline{\J_4}$. Since the
tangent cone has rank two, it is reducible, and thus $\tau$ is a
reducible point.
\end{proof}
Note that in trying to approach the problem in higher genus one
trouble with this corollary is that the locus of reducible ppavs
within $\A_g$ is not contained in the closure of the Jacobian locus,
while for genus 4 this is the case, as we have
$\overline{\J_g}=\A_g$ for $g=1,2,3$.

\section*{Acknowledgements}
We are grateful to Hershel Farkas for drawing our attention to the
problem and encouraging us to work on it. We would also like to
thank Enrico Arbarello and Edoardo Sernesi for discussions on
generalizations to higher genera.

\end{document}